\documentclass[12pt]{article}
\usepackage{verbatim}
\usepackage{fancyhdr}
\usepackage{graphicx}
\usepackage{mathrsfs}
\usepackage{amssymb}
\usepackage{pst-poly}
\usepackage{cite}

\newtheorem{lem}{Lemma}[section]
\newtheorem{thm}{Theorem}[section]

\newenvironment{pf}[1][Proof]{\noindent\textbf{#1.} }{\hfill\rule{1mm}{2mm}}

\makeatletter \@addtoreset{equation}{section} \makeatother

\begin{document}

\title{Roman Bondage Numbers of Some Graphs\thanks {The work was supported by NNSF
of China (No. 11071233).}}
\author
{Fu-Tao Hu, Jun-Ming Xu\footnote{Corresponding
author:
xujm@ustc.edu.cn}\\ \\
{\small Department of Mathematics}  \\
{\small University of Science and Technology of China}\\
{\small Hefei, Anhui, 230026, China} }
\date{}
 \maketitle

\begin{quotation}

\textbf{Abstract}: A Roman dominating function on a graph $G=(V,E)$
is a function $f: V\to \{0,1,2\}$ satisfying the condition that
every vertex $u$ with $f(u)=0$ is adjacent to at least one vertex
$v$ with $f(v)=2$. The weight of a Roman dominating function is the
value $f(G)=\sum_{u\in V} f(u)$. The Roman domination number of $G$
is the minimum weight of a Roman dominating function on $G$. The
Roman bondage number of a nonempty graph $G$ is the minimum number
of edges whose removal results in a graph with the Roman domination
number larger than that of $G$. This paper determines the exact
value of the Roman bondage numbers of two classes of graphs,
complete $t$-partite graphs and $(n-3)$-regular graphs with order
$n$ for any $n\ge 5$.



\end{quotation}

\section{Introduction}

In this paper, a graph $G=(V,E)$ is considered as an undirected
graph without loops and multi-edges, where $V=V(G)$ is the vertex
set and $E=E(G)$ is the edge set. For each vertex $x\in V(G)$, let
$N_G(x)=\{y\in V(G): xy\in E(G)\}$, $N_G[x]=N_G(x)\cup \{x\}$, and
$E_G(x)=\{xy: y\in N_G(x)\}$. The cardinality $|E_G(x)|$ is the
degree of $x$, denoted by $d_G(x)$. For two disjoint nonempty and
proper subsets $S$ and $T$ in $V(G)$, we use $E_G(S,T)$ to denote
the set of edges between $S$ and $T$ in $G$, and $G[S]$ to denote a
subgraph of $G$ induced by $S$.

A subset $D\subseteq V$ is a {\it dominating set} of $G$ if
$N_G(x)\cap D\ne \emptyset$ for every vertex $x$ in $G-D$. The {\it
domination number} of $G$, denoted by $\gamma(G)$, is the minimum
cardinality of all dominating sets of $G$. The {\it Roman dominating
function} on $G$, proposed by Cockayne {\it et al.}~\cite{cd04}, is
a function $f: V\rightarrow\{0,1,2\}$ such that each vertex $x$ with
$f(x)=0$ is adjacent to at least one vertex $y$ with $f(y)=2$. For
$S\subseteq V$ let $f(S)=\sum\limits_{u\in S} f(u)$. The value
$f(V(G))$ is called the {\it weight} of $f$, denoted by $f(G)$. The
{\it Roman domination number}, denoted by $\gamma_{\rm R}(G)$, is
defined as the minimum weight of all Roman dominating functions,
that is,
 $$
 \gamma_{\rm R}(G)=\min \{f(G): f\ {\rm is\ a\ Roman\ dominating\ function\ on}\ G \}.
 $$
A Roman dominating function $f$ is called a {\it $\gamma_{\rm
R}$-function} if $f(G)=\gamma_{\rm R}(G)$. Roman domination numbers
have been studied, see, for example \cite{cd04, ck10, fk09, fy09,
hx10, lk08, lkl05, pp02, sh07, sh10, xc06}.

To measure the vulnerability or the stability of the domination in
an interconnection network under edge failure, Fink et
at.~\cite{fjkr90} proposed the concept of the bondage number in
1990. The {\it bondage number}, denoted by $b(G)$, of $G$ is the
minimum number of edges whose removal from $G$ results in a graph
with larger domination number of $G$.

The {\it Roman bondage number}, denoted by $b_{\rm R}(G)$ and
proposed first by Rad and Volkmann~\cite{rv94}, of a nonempty graph
$G$ is the minimum number of edges whose removal from $G$ results in
a graph with larger Roman domination number. Precisely speaking, the
Roman bondage number
 $$
 b_{\rm R}(G)=\min\{|B|: B\subseteq E(G), \gamma_{\rm R}(G-B)>\gamma_{\rm R}(G)\}.
 $$

An edge set $B$ that $\gamma_{\rm R}(G-B)>\gamma_{\rm R}(G)$ is
called the {\it Roman bondage set} and the minimum one the {\it
minimum Roman bondage set}. In~\cite{hx10}, the authors showed that
the decision problem for $b_{\rm R}(G)$ is NP-hard even for
bipartite graphs.

For a complete $t$-partite graph $K_{m_{1},m_{2},\ldots,m_t}$, its
bondage number determined by Fink {\it et al.}~\cite{fjkr90} for the
undirected case and by Zhang {\it {\it et al.}}~\cite{zlm09} for the
directed case. Motivated by these results, we, in this paper,
consider its Roman bondage number. For a complete $t$-partite
undirected graph $K_{m_{1},m_{2},\ldots,m_t}$ with
$m_1=m_2=\ldots=m_i<m_{i+1}\leq \ldots \leq m_t$ and
$n=\sum\limits_{j=1}^{t} m_j$, we determine that
 $$
 b_{\rm R}(K_{m_{1},m_{2},\ldots,m_t})=\left\{
 \begin{array}{ll}
 \lfloor\frac{i}{2}\rfloor & {\rm if}\ m_i=1~{\rm and}~n\ge 3;\\
 2 \ & {\rm if}\ m_i=2\ {\rm and}\ i=1;\\
 i\  & {\rm if}\ m_i=2\ {\rm and}\ i\geq 2;\\
 n-1 \ & {\rm if}\ m_i=3\ {\rm and}\ i=t\geq 3;\\
 n-m_t & {\rm if}\ m_i\ge 3~{\rm and}~m_t\ge 4.
 \end{array}\right.
 $$
Consider $K_{3,3,\ldots,3}$ of order $n\ge 9$, which is an
$(n-3)$-regular graph. The above result means that $b_{\rm
R}(K_{3,3,\ldots,3})=n-1$. In this paper, we further determine that
$b_{\rm R}(G)=n-2$ for any $(n-3)$-regular graph $G$ of order $n\ge
5$ and $G\ne K_{3,3,\ldots,3}$.

In the proofs of our results, when a Roman dominating function of a
graph is constructed, we only give its nonzero value of some
vertices.

For terminology and notation on graph theory not given here, the
reader is referred to Xu~\cite{x03}.

\section{Preliminary results}

\begin{lem}\label{lem2.1}
{\rm (Cockayne et al. \cite{cd04})}
For a complete $t$-partite graph $K_{m_{1},m_{2},\ldots,m_t}$ with
$1\leq m_1\leq m_2\leq \ldots \leq m_t$ and $t\geq 2$,
 $$
 \gamma_{\rm R}(K_{m_{1},m_{2},\ldots,m_t})=\left\{
 \begin{array}{ll}
 2, & {\rm if}\ m_1=1;\\
 3, & {\rm if}\ m_1=2;\\
 4, & {\rm if}\ m_1\geq 3.
 \end{array}\right.
 $$
\end{lem}

\begin{lem}\label{lem2.2}
{\rm (Hu and Xu \cite{hx10})}
Let $G$ be a graph of order $n\geq 3$ and
$t$ be the number of vertices of degree $n-1$ in $G$. If $t\geq 1$,
then $b_{\rm R}(G)=\lceil\frac{t}{2}\rceil$.
\end{lem}

\begin{lem}\label{lem2.3}
{\rm (Hu and Xu \cite{hx10})}
Let $G$ be a nonempty graph of order $n\geq 3$, then
$\gamma_{\rm R}(G)=3$ if and only if $\Delta(G)=n-2$.
\end{lem}

\begin{lem}\label{lem2.4}
Let $G$ be an {\rm ($n-3$)}-regular graph of order $n\ge 4$. Then
$\gamma_{\rm R}(G)=4$.
\end{lem}
\begin{pf}
Since $G$ is an ($n-3$)-regular graph and $n\ge 4$, $G$ is nonempty.
Let $f$ be a minimum Roman dominating function of $G$. If there is
no vertex $x$ such that $f(x)=2$, then $f(y)=1$ for every vertex
$y$. Therefore, $\gamma_{\rm R}(G)=f(G)\ge 4$. Assume that there is
some vertex $x$ with $f(x)=2$. Let $u$ and $v$ be the only two
vertices not adjacent to $x$ in $G$. If $f(u)=0$ or $f(v)=0$, then
there exists a vertex $y\ne x$ adjacent to $u$ or $v$ in $G$ such
that $f(y)=2$ and hence $\gamma_{\rm R}(G)=f(G)\ge f(x)+f(y)=4$. If
$f(u)\ge 1$ and $f(v)\ge 1$, then $\gamma_{\rm R}(G)=f(G)\ge
f(x)+f(u)+f(v)\ge 4$. In the following, we prove $\gamma_{\rm
R}(G)\le 4$.

For any vertex $x$, let $y$ and $z$ be the only two vertices not adjacent
to $x$ in $G$. Denote $f(x)=2$ and $f(y)=f(z)=1$. Then, $f$ is a Roman
dominating function of $G$ with $f(G)=4$, and hence $\gamma_{\rm R}(G)\leq 4$.
Thus, $\gamma_{\rm R}(G)=4$.
\end{pf}

\begin{lem}\label{lem2.5}
Let $G$ be an {\rm ($n-3$)}-regular graph of order $n\ge 5$ and $B$ be a
Roman bondage set of $G$. Then $E_G(x)\cap B\ne\emptyset$ for any $x\in V(G)$.
\end{lem}

\begin{pf}
By Lemma~\ref{lem2.4}, $\gamma_{\rm R}(G)=4$. Let $G'=G-B$. Then
$\gamma_{\rm R}(G')>4$ since $B$ is a Roman bondage set in $G$. By
contradiction, assume $E_G(x)\cap B=\emptyset$ for some $x\in V(G)$.
Let $y$ and $z$ be the only two vertices not adjacent to $x$ in $G$.
Denote $f(x)=2$ and $f(y)=f(z)=1$. Since every $u\notin \{x,y,z\}$
is adjacent to $x$, $f$ is a Roman dominating function of $G'$ with
$f(G')=4$. Thus we obtain a contradiction as follows. $\gamma_{\rm
R}(G')\leq f(G')=4< \gamma_{\rm R}(G')$.
\end{pf}

\begin{lem}\label{lem2.6}
Let $G$ be an {\rm ($n-3$)}-regular graph of order $n\ge 5$ and $B$
be a Roman bondage set of $G$, $x$ be any vertex,  $y$ and $z$ be
the only two vertices not adjacent to $x$ in $G$. If $E_G(x)\cap
 B=\{xw\}$, then $|E_G(\{y,z,w\},x')\cap B|\ge 1$ for any vertex
$x'\in V(G)\setminus \{x,y,z,w\}$ that is adjacent to each vertex in
$\{y,z,w\}$ in $G$.
\end{lem}

\begin{pf}
Since $E_G(x)\cap B\ne\emptyset$ by Lemma~\ref{lem2.5}, there is a
vertex $w$ such that $xw\in E_G(x)\cap B$. Let $G'=G-B$. Then
$\gamma_{\rm R}(G')>4$ by Lemma~\ref{lem2.4}. By contradiction,
suppose $E_G(\{y,z,w\},x')\cap B=\emptyset$ for some vertex $x'\in
V(G)\setminus \{x,y,z,w\}$ that is adjacent to each in $\{y,z,w\}$
in $G$. Set $f(x)=f(x')=2$. Then, $f$ is a Roman dominating function
of $G'$ with $f(G')=4$ since $N_{G'}[x]\cup N_{G'}[x']=V(G)$, a
contradiction.
\end{pf}

\begin{lem}\label{lem2.7}
Let $G$ be an {\rm ($n-3$)}-regular graph of order $n\ge 6$ and $B$
be a Roman bondage set of $G$. For three vertices $x,y$ and $z$ that
are pairwise non-adjacent in $G$, if each of them is incident
with exact one edge in $B$, then $|B|\geq n-2$ and, moreover,
$|B|\ge n-1$ if $G=K_{3,3,\ldots,3}$.
\end{lem}

\begin{pf}
By the hypothesis, for any $v\in\{x,y,z\}$, $|E_G(v)\cap B|=1$ and
$v$ is adjacent to every vertex in $V(G-\{x,y,z,v\})$. Let $xu\in
E_G(x)\cap B$. We claim $yu\in E_G(y)\cap B$ and $zu\in E_G(z)\cap
B$. In fact, by contradiction, without loss of generality suppose
$yv\in E_G(y)\cap B$ and $zw\in E_G(z)\cap B$ with $u\ne v$ and
$u\ne w$. Then, $u$ is adjacent to $y$ and $z$ in $G-B$. Set
$f(x)=f(u)=2$. Then $f$ is a Roman dominating function of $G$ with
$f(G-B)=4$, which contradicts with $\gamma_R(G-B)>4$ by
Lemma~\ref{lem2.4}.

Let $s$ and $t$ be the only two vertices not adjacent to $u$ in $G$,
and let $V'=V(G)\setminus \{x,y,z,u,s,t\}$. By the hypothesis, each
vertex in $\{y,z,u\}$ is adjacent to all vertices in $V'$ in $G$. By
Lemma~\ref{lem2.6}, for any vertex $x'\in V'$, if such a vertex
exists, $|E_G(\{u,y,z\},x')\cap B|\ge 1$, and so
 \begin{equation}\label{e2.1}
 |E_G(\{u,y,z\},V')\cap B|\ge |V'|=n-6.
 \end{equation}
By Lemma~\ref{lem2.5}, $|E_G(s)\cap B|\ge 1$ and $|E_G(t)\cap B|\ge
1$, and so we have that
 \begin{equation}\label{e2.2}
  |(E_G(s)\cup E_G(t))\cap B|\ \left\{\begin{array}{ll}
  \ge 1\ \ & {\rm if}\ \ st\in E(G);\\
  =2\ \ & {\rm if}\ \ st\notin E(G).
  \end{array}\right.
  \end{equation}
 It follows from (\ref{e2.1}) and (\ref{e2.2}) that
 $$
 \begin{array}{rl}
 |B| &\ge |\{xu,yu,zu\}|+|(E_G(s)\cup E_G(t))\cap  B|\\
     &\qquad +|E_G(\{u,y,z\},V')\cap B|\\
     &\ge \left\{\begin{array}{ll}
      n-2\ \ & {\rm if}\ \ st\in E(G);\\
      n-1\ \ & {\rm if}\ \ st\notin E(G).
  \end{array}\right.
     \end{array}
 $$
If $G=K_{3,3,\ldots,3}$, then $|(E_G(s)\cup E_G(t))\cap B|\ge 2$
since then $st\notin E(G)$ and, hence, $|B|\ge n-1$.
\end{pf}

\begin{lem}\label{lem2.8}
Let $G$ be an {\rm ($n-3$)}-regular graph of order $n\ge 5$ and $B$
be a Roman bondage set of $G$. Let $x\in V(G)$, $y$ and $z$ be the
only two vertices not adjacent to $x$ in $G$. Let $E_G(x)\cap
B=\{xw\}$ and $G'=G-B$. Then $|E(G'[\{y,z,w\}])|\le 1$, that is,
 $$
 |E(G[\{y,z,w\}])\cap B|\ge \left\{\begin{array}{ll}
  1\ \ & {\rm if}\ \ |E(G[\{y,z,w\}])|=2;\\
  2\ \ & {\rm if}\ \ |E(G[\{y,z,w\}])|=3.
  \end{array}\right.
 $$
\end{lem}

\begin{pf}
Suppose to the contrary that $|E(G'[\{y,z,w\}])|\ge 2$. Without loss
of generality, let $yw,zw\in E(G')$. Denote $f(x)=f(w)=2$. Note that
$x$ is adjacent to every vertex except $w$, $y$ and $z$ in $G'$.
Thus, $f$ is a Roman dominating function of $G'$ with $f(G')=4$, a
contradiction with $\gamma_{\rm R}(G')>4$ by Lemma~\ref{lem2.4}.
\end{pf}

\section{Results on complete $t$-partite graphs}

\begin{thm}\label{thm3.1}\
Let $G = K_{m_1,m_2,\ldots,m_t}$ be a complete $t$-partite graph
with $m_1=m_2=\ldots=m_i<m_{i+1}\leq \ldots \leq m_t$, $t\geq 2$
and $n=\sum\limits_{j=1}^{t} m_j$. Then
 $$
 b_{\rm R}(G)=\left\{\begin{array}{ll}
 \lceil\frac{i}{2}\rceil\  & {\rm if}\ m_i=1\  {\rm and}\ n\geq 3;\\
 2 \ & {\rm if}\ m_i=2\ {\rm and}\ i=1;\\
 i\  & {\rm if}\ m_i=2\ {\rm and}\ i\geq 2;\\
 n-1 \ & {\rm if}\ m_i=3\ {\rm and}\ i=t\geq 3;\\
 n-m_t  \ & {\rm if}\ m_i\geq 3\  {\rm and}\ m_t\geq 4.
 \end{array}\right.
 $$
\end{thm}

\begin{pf}
Let $\{X_1,X_2,\ldots, X_t\}$ be the corresponding $t$-partitions of
$V(G)$.

(1) If $m_i=1$ and $n\geq 3$, then $G$ has $i$ vertices of degree
$n-1$, and so $b_{\rm R}(G)=\lceil\frac{i}{2}\rceil$ by
Lemma~\ref{lem2.2}.

\vskip6pt

(2) If $m_i=2$, then $\Delta(G)=n-2$. By Lemma~\ref{lem2.1},
$\gamma_{\rm R}(G)=3$. Let $B\subseteq E(G)$ be a Roman bondage set
of $G$ with $|B|= b_{\rm R}(G)$ and $G'=G-B$. Then $\gamma_{\rm
R}(G')>\gamma_{\rm R}(G)=3$, and so $\Delta(G')\leq n-3$ by
Lemma~\ref{lem2.3}. Thus, $|B\cap E_G(x)|\geq 1$ for every vertex in
$X_j$ $(1\leq j\leq i)$, that is, $|B|\geq 2$ if $i=1$ and $|B|\geq
i$ if $i>1$.

If $i=1$, then only two vertices of degree $n-2$ are in $X_1$, and
the removal of any two edges incident with distinct vertices in
$X_1$ results in a graph $G''$ with $\Delta(G'')\leq n-3$, and hence
$\gamma_{\rm R}(G'')\neq 3$ by Lemma~\ref{lem2.3}. Since
$\gamma_{\rm R}(G'')\geq \gamma_{\rm R}(G)=3$, $\gamma_{\rm
R}(G'')\geq 4$. Thus, $b_{\rm R}(G)\leq 2$ and hence $b_{\rm
R}(G)=2$.

If $i>1$, then the subgraph $H$ induced by
$\bigcup\limits_{j=1}^{i}{X_j}$ of $G$ is a complete $i$-partite
graph with each partition consisting of two vertices, which is
$2$-edge-connected and $2(i-1)$-regular, and so has a perfect
matching $M$ with $|M|=i$. Thus, $G-M$ has the maximum degree $n-3$
and hence $\gamma_{\rm R}(G-M)\neq 3$ by Lemma~\ref{lem2.3}. Since
$\gamma_{\rm R}(G-M)\geq \gamma_{\rm R}(G)=3$, $\gamma_{\rm
R}(G-M)\geq 4$. Thus, $b_{\rm R}(G)\leq |M|=i$, and so $b_{\rm
R}(G)=i$.

\vskip6pt

(3) Assume $m_i=3$ and $i=t$. Then $G$ is $(n-3)$-regular. Note that
if $t=2$ then $n=6$ and $b_R(K_{3,3})=4=n-2$. It is easy to verify
that the conclusion is true for $t=3,4$, and assume $t\geq 5$ below.
Let $x\in V(G)$ and $H=G-E_G(x)$, then $\gamma_{\rm
R}(H)=1+\gamma_{R}(K_{2,3,\ldots,3})=4$ by Lemma~\ref{lem2.1}. By
the conclusion (2) just showed, $b_{\rm R}(K_{2,3,\ldots,3})=2$ and
hence
 $$
 b_{\rm R}(G)\leq |E_G(x)|+b_{\rm R}(K_{2,3,\ldots,3})=(n-3)+2=n-1.
 $$

We now prove that $b_{\rm R}(G)\geq n-1$. By contradiction. Assume
that there is a Roman bondage set $B$ of $G$ such that $|B|\le n-2$.
Let $G'=G-B$. Then $\gamma_{\rm R}(G')>\gamma_{\rm R}(G)=4$ by
Lemma~\ref{lem2.1}, and $|E_G(x)\cap B|\geq 1$ for any vertex $x\in
V(G)$ by Lemma~\ref{lem2.5}. If $|E_G(x)\cap B|\geq 2$ for any
vertex $x\in V(G)$, then the subgraph induced by $B$ has the minimum
degree at least two, and so $|B|\geq n$, a contradiction. Thus,
there is a vertex $x_1$ in $G$ such that $|E_G(x_1)\cap B|=1$. Let
$x_1y_1\in B$ and, without loss of generality, let
$X_1=\{x_1,x_2,x_3\}$ and $X_2=\{y_1,y_2,y_3\}$.
By Lemma~\ref{lem2.8},
 \begin{equation}\label{e3.1}
 |E(G[\{y_1,x_2,x_3\}])\cap B|\ge 1,
 \end{equation}
and by Lemma~\ref{lem2.5},
 \begin{equation}\label{e3.2}
 |E_G(y_2)\cap B|\geq 1\ \ {\rm and}\ |E_G(y_3)\cap B|\geq 1.
 \end{equation}
Let $V_1=V(G)\setminus (X_1\cup X_2)$. By Lemma~\ref{lem2.6},
 \begin{equation}\label{e3.3}
 |E_G(\{y_1,x_2,x_3\},x')\cap B|\ge 1\ \ {\rm for\ any} \ x'\in V_1,
 \end{equation}
and so
 \begin{equation}\label{e3.4}
 |E_G(\{y_1,x_2,x_3\},V_1)\cap B|\ge n-6.
 \end{equation}
It follows from (\ref{e3.1}), (\ref{e3.2}) and (\ref{e3.4}) that
 \begin{equation}\label{e3.5}
 \begin{array}{rl}
 n-2\geq |B|&\geq |\{x_1y_1\}|+|E(G[\{y_1,x_2,x_3\}])\cap B|\\
  &\qquad +|E_G(\{y_1,x_2,x_3\},V_1)\cap B|+|E_G(y_2)\cap B|\\
  &\qquad +|E_G(y_3)\cap B|+|E(G[V_1])\cap B|\\
  &\geq 1+1+(n-6)+1+1+0\\
  &\geq n-2.
  \end{array}
 \end{equation}
Thus, all the equalities in (\ref{e3.5}) hold, which implies that
all the equalities in (\ref{e3.1}), (\ref{e3.2}) and (\ref{e3.3})
hold, and $|E(G[V_1])\cap B|=0$.

Let $E_G(y_2)\cap B=\{y_2u\}$ and $E_G(y_3)\cap B=\{y_2v\}$. The
worst case is that $u$ and $v$ belong to different partitions of
$X_3,\ldots,X_t$. Since $t\ge 5$, there exists some $i$ with $3\le
i\le t$ such that both $u$ and $v$ are not belong to $X_i$. Thus,
each vertex in $X_i$ is incident with exact one edge in $B$. By
Lemma~\ref{lem2.7}, $|B|\ge n-1$, a contradiction.

Thus,  $b_{\rm R}(K_{3,3,\ldots,3})=n-1$.

\vskip6pt

(4) We now assume $m_i\geq 3$ and $m_t\geq 4$. By
Lemma~\ref{lem2.1}, we have $\gamma_{\rm R}(G)=4$. Let $u$ be a
vertex in $X_t$ and $f$ be a $\gamma_{\rm R}$-function of
$G-E_G(u)$. Then $u$ is an isolated vertex. Thus, $f(u)=1$ and
$f(G-u)=4$ by Lemma~\ref{lem2.1} since $G-u$ is a complete
$t$-partite graph with at least 3 vertices in every partition. Thus
$\gamma_{\rm R}(G-E_G(u))=5>4= \gamma_{\rm R}(G)$, and hence $b_{\rm
R}(G)\leq |E_G(u)|=n-m_t$.

We now show $b_{\rm R}(G)\geq n-m_t$. Let $B$ be a minimum Roman
bondage set of $G$, and let $G'=G-B$. Then $\gamma_{\rm
R}(G')>\gamma_{\rm R}(G)=4$.

Assume that there is a vertex $x$ in $G$ such that $E_G(x)\cap
B=\emptyset$. Then there is some $j$ with $1\leq j\le t$ such that
$x\in X_j$. If there exists some $y\in V(G-X_j)$ such that
$E_G(y,X_j)\cap B=\emptyset$. Set $f(x)=f(y)=2$. Then $f$ is a Roman
dominating function of $G'$ with $f(G')=4$, a contradiction. Thus,
 $$
 E_G(y,X_j)\cap B\neq \emptyset\ \ {\rm for\ any} \ y\in V(G-X_j).
 $$
It follows that
 $$
 |B|\geq |V(G)\setminus X_j|=n-m_j\geq n-m_t.
 $$

We now assume that
 \begin{equation}\label{e3.6}
 |E_G(x)\cap B|\geq 1\ \ {\rm for\ any} \ x\in V(G).
 \end{equation}

If $|E_G(x)\cap B|\geq 2$ for any  $x\in V(G)$, then the subgraph
induced by $B$ has the minimum degree at least two, from which we
have $|B|\geq n>n-m_t$.

We suppose that there exists a vertex $x_1\in V(G)$ such that
$|E_G(x_1)\cap B|=1$. Let $x_1\in X_j$ and $x_2,x_3,\ldots,x_{m_j}$
be the other vertices of $X_j$. Let $y_1$ be the unique neighbor of
$x_1$ in $E_G(x_1)\cap B$, and let $X_k$ contains $y_1$. Let
$V'=V(G)\setminus (X_j\cup X_k)$ and $V''=\{y_1,x_2,
x_3,\ldots,x_{m_j}\}$. If there is some $x'\in V'$ such that
$|E_G(x',V'')\cap B|=0$, set $f(x)=f(x')=2$, then $f$ is a Roman
dominating function of $G'$ with $f(G')=4$, a contradiction. Thus,
\begin{equation}\label{e3.7}
 |E_G(x',V'')\cap B|\geq 1\ \ {\rm for\ any} \ x'\in V'.
 \end{equation}
It follows from (\ref{e3.6}) and (\ref{e3.7}) that
 $$
 b_{\rm R}(G)=|B|\geq |V'|+|X_k|\geq n-m_t.
 $$
Thus, $b_{\rm R}(G)=n-m_t$.

The theorem follows.
\end{pf}



\section{Results on $(n-3)$-regular graphs}

By Theorem~\ref{thm3.1}, we immediately have $b_{\rm
R}(K_{3,3,\ldots,3})=n-1$ if its order is $n$. $K_{3,3,\ldots,3}$ is
an $(n-3)$-regular graph of order $n$. In this section, we determine
the Roman bondage number of any $(n-3)$-regular graph $G$ of order
$n$ is equal to $n-2$ if $G\ne K_{3,3,\ldots,3}$.

\begin{lem}\label{lem4.1}
Let $G$ be an {\rm ($n-3$)}-regular graph of order $n\ge 7$ but $G\ne K_{3,3,\ldots,3}$
and $B$ be a Roman bondage set of $G$. Let $x,w\in V(G)$ and $xw\in E(G)$. Let $y,z$
and $p,q$ be the only two vertices not adjacent to $x$ and $w$ in $G$, respectively.
If $E_G(x)\cap B=\{xw\}$ and $\{y,z\}\cap \{p,q\}\ne \emptyset$, then $|B|\geq n-2$.
\end{lem}
\begin{pf}
By Lemma~\ref{lem2.4}, $\gamma_{\rm R}(G)=4$. Let $G'=G-B$. Then $\gamma_{\rm R}(G')>4$.
By Lemma~\ref{lem2.5}, $E_G(y')\cap B\ne\emptyset$ for any $y'\in V(G)$. By contradiction,
assume $|B|\le n-3$. We now deduce a contradiction by considering the following two
cases.

\begin{description}

\item [Case 1]  $\{y,z\}=\{p,q\}$.

In this case, $yz\in E(G)$ since $G$ is $(n-3)$-regular. Let
$U_1=V(G)\setminus \{x,y,z,w\}$. Then any vertex in $U_1$ is
adjacent to each in $\{w,y,z\}$. By Lemma~\ref{lem2.6}, we have that
$|E_G(\{w,y,z\},x')\cap B|\ge 1$, and so $|E_G(\{w,y,z\},U_1)\cap
B|\geq |U_1|=n-4$. It follows that
 \begin{equation}\label{e4.1}
 \begin{array}{rl}
 n-3\geq |B|&\geq |\{xw\}|+|E_G(\{w,y,z\},U_1)\cap B|+|E(G[U_1])\cap B|\\
            &\geq 1+(n-4)+0\\
            &=n-3.
  \end{array}
 \end{equation}
This means that all equalities in (\ref{e4.1}) hold, that is,
$yz\notin B$, $E(G[U_1])\cap B=\emptyset $, $|E_G(\{w,y,z\}, x')\cap
B|=1$ and then, $|E_G(x')\cap B|=1$ for any vertex $x'\in U_1$. Let
$yr\in B$ for some $r\in U_1$ since $E_G(y)\cap B\ne \emptyset$, $s$
and $t$ be the only two vertices not adjacent to $r$ in $G$.

Assume $st\notin E(G)$. Then $r,s,t$ be three vertices not adjacent
to each other in $G$, and each one of them is incident with exact one
edge in $B$. By Lemma~\ref{lem2.7}, $|B|\ge n-2$, a contradiction.

Now, assume $st\in E(G)$. We claim that $ys\in B$ and $yt\in B$.
By contradiction, assume $ys\notin B$. Denote $f(r)=f(s)=2$. Then,
$f$ is a Roman dominating function of $G'$ with $f(G')=4$, a contradiction.
Also, $yt\in B$ by replace $t$ by $s$. Then $zs$ and $zt$ not belong to $B$.
Denote $f(r)=f(z)=2$. Then, $f$ is a Roman dominating function of $G'$ with
$f(G')=4$, a contradiction.

\item [Case 2]  $|\{y,z\}\cap \{p,q\}|=1$. Without loss of generality, let $p=y$.

In this case, $yz,wz\in E(G)$ and hence $|E(G[\{y,z,w\}])\cap B|\geq
1$ by Lemma~\ref{lem2.8}. Let $r$ be the only vertex except $x$ not
adjacent to $z$ in $G$. By Lemma~\ref{lem2.6},
$|E_G(\{w,y,z\},x')\cap B|\geq 1$ for any vertex $x'\in
U_2=V(G)\setminus \{x,y,z,w,q,r\}$.

If $q=r$, then $|E_G(\{w,y,z\},U_2)\cap B|\geq |U_2|=n-5$. Then we
can deduce a contradiction as follows.

$$
 \begin{array}{rl}
 n-3\geq |B|&\geq |\{xw\}|+|E_G(\{w,y,z\},U_2)\cap B|\\
            &\qquad+E(G[\{y,z,w\}])\cap B|+|E_G(q)\cap B|\\
            &\geq 1+(n-5)+1+1\\
            &=n-2.
  \end{array}
$$

If $q\neq r$, then $wr,zq\in E(G)$ and $|E_G(\{w,y,z\},U_2)\cap
B|\geq |U_2|=n-6$. Then,

\begin{equation}\label{e4.2}
 \begin{array}{rl}
 n-3\geq |B|&\geq |\{xw\}|+|E_G(\{w,y,z\},U_2)\cap B|+|E(G[U_2])\cap B|\\
            &\qquad+E(G[\{y,z,w\}])\cap B|+|(E_G(q)\cup E_G(r))\cap B|\\
            &\geq 1+(n-6)+0+1+1\\
            &=n-3.
  \end{array}
\end{equation}

It follows that the equalities in (\ref{e4.2}) hold, which implies
that $|(E_G(q)\cup E_G(r))\cap B|=1$, $E(G[U_2])\cap B=\emptyset$,
$|E_G(\{w,y,z\},x')\cap B|=1$ and then, $|E_G(x')\cap B|=1$ for any
vertex $x'\in U_2$. Then $(E_G(q)\cup E_G(r))\cap B=\{qr\}$, and
hence $wr\notin B$, $zq\notin B$.

Let $s$ be the only vertex except $w$ not adjacent to $q$ in $G$.
Then both $rs$ and $ws$ not belong to $G'$, otherwise denote $f(q)=f(r)=2$ or
$f(q)=f(w)=2$. Then $f$ is a Roman dominating function of $G'$ with $f(G')=4$,
a contradiction. $rs, ws\notin E(G')$ imply that $ws\in B$ and $rs\notin E(G)$.
Then $zs\in E(G)$ and $zs\notin B$ since $|E_G(\{w,y,z\},s)\cap B|=1$. Denote
$f(r)=f(z)=2$. Then $f$ is a Roman dominating function of $G'$ with $f(G')=4$,
a contradiction. Thus, $|B|\geq n-2$.

\end{description}
The lemma follows. \end{pf}

\begin{lem}\label{lem4.2}
let $G$ be an {\rm ($n-3$)}-regular graph of order $n\ge 7$
but $G\ne K_{3,3,\ldots,3}$ and $B$ be a Roman bondage set of $G$.
Let $x,w\in V(G)$ and $xw\in E(G)$. If $E_G(x)\cap B= E_G(w)\cap B=\{xw\}$,
then $|B|\geq n-2$.
\end{lem}

\begin{pf}
Let $y,z$ and $p,q$ be the only two vertices not adjacent to $x$ and
$w$ in $G$, respectively.

We claim that $\{y,z\}\cap \{p,q\}\ne \emptyset$. By contradiction,
suppose $\{y,z\}\cap \{p,q\}=\emptyset$. Then $wy,wz\in E(G)$, and
$wy,wz\notin B$ since $E_G(w)\cap B=\{xw\}$. Denote $f(x)=f(w)=2$.
Then $f$ is a Roman dominating function of $G'$ with $f(G')=4$,
a contradiction. Thus $\{y,z\}\cap \{p,q\}\ne \emptyset$, and hence
$|B|\geq n-2$ by Lemma~\ref{lem4.1}.
\end{pf}

\begin{thm}\label{thm4.1}
Let $G$ be an {\rm ($n-3$)}-regular graph of order $n \ge 5$  but $G\ne K_{3,3,\ldots,3}$.
Then $b_{\rm R}(G)=n-2$.
\end{thm}
\begin{pf}
We first consider $n\in\{5,6\}$. If $n=5$, then $G=C_5$, and so
$b_{\rm R}(G)=3$. If $n=6$, then $G$ is the Cartesian product of a
cycle $C_3$ and a complete graph $K_2$, that is, $G=C_3\times K_2$,
and so $b_{\rm R}(G)=4$. In the following, suppose $n\ge 7$.

By Lemma~\ref{lem2.4}, $\gamma_{\rm R}(G)=4$. Let  $x_0\in V(G)$ and
$y_0z_0\in E(G)$, where $y_0$ and $z_0$ are the only two vertices
not adjacent to $x_0$ in $G$. We consider the Roman domination
number of $H=G-x_0-y_0z_0$. Since $H$ is $(|V(H)|-3)$-regular and
$|V(H)|\ge 4$, $\gamma_{\rm R}(H)=4$ by Lemma~\ref{lem2.4}. Thus
$\gamma_{\rm R}(G-E_G(x_0)-y_0z_0)\ge 5$ and hence $b_{\rm R}(G)\leq
|E_G(x_0)|+1=n-2$. Next, we prove that $b_{\rm R}(G)\geq n-2$.

Let $B$ be a minimum Roman bondage set of $G$ and $G'=G-B$. Then
$|B|\leq n-2$ and $\gamma_{\rm R}(G')>4$. We now prove $|B|\ge n-2$.
By contradiction, assume $|B|\le n-3$. By Lemma~\ref{lem2.5},
$E_G(y')\cap B\ne\emptyset$ for any $y'\in V(G)$. Then there exists
a vertex $x$ such that $|E_G(x)\cap B|=1$. Let $xw\in B$, $y$ and
$z$ be the only two vertices not adjacent to $x$ in $G$. Let $p$ and
$q$ be the only two vertices not adjacent to $w$ in $G$. If
$\{y,z\}\cap \{p,q\}\ne \emptyset$, then $|B|\geq n-2$ by
Lemma~\ref{lem4.1}. Thus, we only need to consider the case of
$\{y,z\}\cap \{p,q\}=\emptyset$. In this case, $wy,wz\in E(G)$. We
now deduce a contradiction by considering the following two cases.

\begin{description}

\item [Case 1] $yz\notin E(G)$.

By Lemma~\ref{lem2.8}, $|E(G[\{y,z,w\}])\cap B|\ge 1$. By Lemma~\ref{lem2.6},
$|E_G(\{w,y,z\},x')\cap B|\geq 1$ for any vertex $x'\in X_1=V(G)\setminus \{x,y,z,w,p,q\}$,
and so $|E_G(\{w,y,z\},X_1)\cap B|\geq |X_1|=n-6$. Then,
\begin{equation}\label{e3.7a}
 \begin{array}{rl}
 n-3\geq |B|&\geq |\{xw\}|+|E_G(\{w,y,z\},X_1)\cap B|\\
            &\qquad+|E(G[\{y,z,w\}])\cap B|+|(E_G(p)\cup E_G(q))\cap B|\\
            &\geq 1+(n-6)+1+1\\
            &=n-3.
  \end{array}
\end{equation}

It follows that the equalities in (\ref{e3.7a}) hold, which implies
that $|E_G(\{p,q\})\cap B|=1$. Then $(E_G(p)\cup E_G(q))\cap
B=\{pq\}$ and then, $E_G(p)\cap B=E_G(q)\cap B=\{pq\}$. By
Lemma~\ref{lem4.2}, $|B|\ge n-2$, a contradiction.

\item [Case 2]  $yz\in E(G)$.

Let $r$ and $s$ be the only vertex except $x$ not adjacent to $y$
and $z$ in $G$, respectively. By Lemma~\ref{lem2.8},
$|E(G[\{w,y,z\}])\cap B|\geq 2$. By Lemma~\ref{lem2.6},
$|E_G(\{w,y,z\},x')\cap B|\geq 1$ for any vertex $x'\in
X_2=V(G)\setminus\{x,y,z,w,p,q,r,s\}$. Thus, we have that
 \begin{equation}\label{e4.4}
 |E_G(\{w,y,z\},X_2)\cap B|\ge |X_2|\ge
 \left\{\begin{array}{ll}
 n-6\ &\ {\rm if}\ |\{r,s\}\cup \{p,q\}|\le 2;\\
 n-7\ &\ {\rm if}\ |\{r,s\}\cup \{p,q\}|=3;\\
 n-8\ &\ {\rm if}\ |\{r,s\}\cup \{p,q\}|=4,
 \end{array}\right.
 \end{equation}
and
 \begin{equation}\label{e4.5}
 |(E_G(p)\cup E_G(q)\cup E_G(r)\cup E_G(s))\cap B|\geq
 \left\{\begin{array}{ll}
 1\ &\ {\rm if}\ |\{r,s\}\cup \{p,q\}|\le 2;\\
 2\ &\ {\rm if}\ |\{r,s\}\cup \{p,q\}|=3;\\
 2\ &\ {\rm if}\ |\{r,s\}\cup \{p,q\}|=4,
 \end{array}\right.
 \end{equation}
It follows from (\ref{e4.4}) and (\ref{e4.5}) that
 \begin{equation}\label{e4.6}
 \begin{array}{rl}
n-3\ge |B|&\geq |\{xw\}|+|E_G(\{w,y,z\},X_2)\cap B|+|E(G[\{w,y,z\}])\cap B|\\
            &\qquad +|(E_G(p)\cup E_G(q)\cup E_G(r)\cup E_G(s))\cap B|\\
           &\geq  \left\{\begin{array}{ll}
            n-2\ \ & {\rm if}\ |\{r,s\}\cup \{p,q\}|\le 3;\\
            n-3\ \ & {\rm if}\ |\{r,s\}\cup \{p,q\}|=4.
            \end{array}\right.
  \end{array}
 \end{equation}

The Eq. (\ref{e4.6}) implies that $|\{r,s\}\cup \{p,q\}|=4$,
$|B|=n-3$ and $|(E_G(p)\cup E_G(q)\cup E_G(r)\cup E_G(s))\cap B|=2$.
Then there exist two vertices $u,v$ in $\{p,q,r,s\}$ such that
$E_G(u)\cap B=E_G(v)\cap B=\{uv\}$. By Lemma~\ref{lem4.2}, $|B|\ge
n-2$, a contradiction.

\end{description}

Thus, $b_{\rm R}(G)=n-2$, and so the theorem follows.
\end{pf}

\end{document}